\documentclass[]{article}
\usepackage{amsmath}
\usepackage{amssymb}
\usepackage{hyperref}

\numberwithin{equation}{section}

\newtheorem{thm}{Theorem}[section]
\newtheorem{cor}[thm]{Corollary}

\newtheorem{lemma}[thm]{Lemma}

\newtheorem{theorem}[thm]{Theorem}

\newenvironment{proofof}{}{\hfill$\square$\vskip.5cm}

\newcommand{\C}{\mathbb{C}}

\newcommand{\N}{\mathbb{N}}

\title{A Universality Property of Gaussian Analytic Functions}

\author{Andrew Ledoan\footnote{Department of Mathematics, University of Rochester, Rochester, NY 14627, USA, [ledoan,~sstarr]@math.rochester.edu.  The work of S.S. was supported in part by a U.S. National Science Foundation grant DMS-0706927.}\, , Marco Merkli\footnote{Department of Mathematics and Statistics, Memorial University of Newfoundland, St. John's, NL, Canada, A1C 5S7, merkli@mun.ca. Suppported by the Natural Sciences and Engineering Research Council of Canada (NSERC) under Grant 205247.} \  and Shannon Starr${}^*$}

\date{August 16, 2010}

\begin{document}

\maketitle

\begin{abstract}
\setcounter{section}{0}

We consider random analytic functions defined on the unit disk of the complex plane as power series such that the coefficients
are i.i.d., complex valued random variables, with mean zero and unit variance.
For the case of complex Gaussian coefficients, Peres and Vir\'ag showed that the zero
set forms a determinantal
point process with the Bergman kernel.
We show that for general choices of random coefficients, the zero set is asymptotically given
by the same distribution near the boundary of the disk, which expresses a universality property.
The proof is elementary and general.

\vspace{8pt}
\vskip .2 cm
\noindent
{\small \bf MCS numbers: 30B20, 60B12, 60G15}
\end{abstract}

\section{Main Result}

Random analytic functions are a topic of classical interest \cite{B-RS, Kahane},
which  gained renewed interest,
as a toy model for quantum chaos
following work of Bogomolny, Bohigas and Leboeuf \cite{BBL1,BBL2}.
A recent short paper about Gaussian analytic functions
is entitled, ``What is \dots a Gaussian entire function,'' \cite{NazarovSodin}.

Given a sequence of coefficients
$\boldsymbol{x} = (x_0,x_1,x_2,\dots)$,
one may define the power series
$$
f(\boldsymbol{x},z)\, =\, \sum_{n=0}^{\infty} x_n z^n\, .
$$
We consider random analytic functions defined by choosing a coefficient
sequence $\boldsymbol{X} = (X_0,X_1,\dots)$ where $X_0, X_1, \dots$ are i.i.d.,
complex valued random variables, with mean zero and unit variance, such that
\begin{equation}
\label{eq:unitvariance}
{\bf E}\big[ (\textrm{Re}[X_i])^2\big]\, 
=\, {\bf E} \big[ (\textrm{Im}[X_i])^2\big]\, ,\quad
{\bf E} \textrm{Re}[X_i] \textrm{Im}[X_i]\, =\, 0\, .
\end{equation}
By the Borel-Cantelli lemma, the radius of convergence of $f(\boldsymbol{X},z)$ is 1,
almost surely.

We wish to consider the case where the $X_i$'s are not necessarily Gaussian.
But we first recall a  beautiful result
of Peres and Vir\'ag for Gaussian analytic functions \cite{PeresVirag}.

\begin{theorem}[Peres and Vir\'ag, 2005]
\label{thm:PV}
Suppose that $X_0,X_1,\dots$ are i.i.d., and each $X_i$
has density on $\C$ given by
$\pi^{-1} \exp(-|z|^2)$.
Then the zero set is a determinantal point process with Bergman kernel,
$$
K(z_1,z_2)\, =\, \frac{\pi}{(1-z_1 \overline{z}_2)^2}\, .
$$
\end{theorem}

This kernel is invariant under the action of the symmetries of the hyperbolic plane
modeled by the Poincar\'e disk.
For each $u \in \C$ with $|u|<1$, the M\"obius transformation is one such isometry
$$
\Phi(u;z)\, =\, \frac{z-u}{1-\overline{u} z}\, .
$$
This maps the open unit disk $U(0,1) = \{z \in \C \, :\, |z|<1 \}$
bijectively onto itself.
Note that for a fixed $z \in U(0,1)$, 
$$
|\Phi(u,z)| \uparrow 1\quad \textrm{as}\quad |u| \uparrow 1\, .
$$
In other words, taking $|u| \uparrow 1$ maps every point in the interior of the disk to a single point on the boundary.
Our main result establishes a limit law at this boundary, for general random coefficients,
not necessarily Gaussian.

\begin{theorem}
\label{thm:main}
Suppose that $X_0,X_1,\dots$ are i.i.d., complex valued random variables with mean zero
and satisfying (\ref{eq:unitvariance}).
Let $Z(\boldsymbol{X})$ denote the random zero set $\{ \xi \in U(0,1)\, :\, f(\boldsymbol{X},\xi)=0\}$.
Then, for any $n \in \N$, and any $n$ distinct points $z_1,\dots,z_n \in U(0,1)$
$$
\lim_{\epsilon \downarrow 0} \lim_{|u| \uparrow 1}
\epsilon^{-2n} {\bf P}\left(\bigcap_{i=1}^{n} \{U(z_i,\epsilon) \cap Z(\boldsymbol{X}) \neq \emptyset \} \right)\,
=\, \det\big( K(z_i,z_j) \big)_{i,j=1}^{n}\, ,
$$
where we write $U(z_i,\epsilon)$ for the open ball $\{z \in \C \, : \, |z-z_i| < \epsilon \}$.
\end{theorem}

There have been many papers proving convergence of the first intensity measure of zeroes.
In fact there are very precise and general results in this direction.
See for example \cite{EdelmanKostlan,IbragimovZeitouni,SheppVanderbei}.
But this result addresses a slightly different issue because in principle it also gives
correlations.

There is another important group of papers proving universality for Gaussian analytic functions,
for the entire ensemble of correlations, by Bleher, Shiffman and Zelditch \cite{BD,BSZ,BSZ2}.
But these are in a different context.

In Section 2 we give the simple proof of this result.
In Section 3 we describe the extensions to a related family of Gaussian analytic functions
considered by Hough, Krishnapur, Peres and Vir\'ag in their recent monograph \cite{HKPV},
which have interesting properties but whose zero sets are not determinantal.

\section{Proof of the Main Result}

A main step in proving Theorem \ref{thm:main} is the following elementary observation.

\begin{lemma}
\label{lem:clt}
Let $\boldsymbol{Y} = (Y_0,Y_1,\dots)$ be i.i.d., complex Gaussians such that each $Y_i$
has density equal to $\pi^{-1} e^{-|y|^2}$ on $\C$.
Then for any $n \in \N$, any $z_1,\dots,z_n \in U(0,1)$ and any $\lambda_1,\dots,\lambda_n \in \C$, we have
$$
\textrm{Re}\left[ \sum_{i=1}^{n} \lambda_i \frac{f(\boldsymbol{X},\Phi(u,z_i))}{\Delta(u,z_i)} \right]\,
\Rightarrow\, \textrm{Re}\left[ \sum_{i=1}^{n} \lambda_i f(\boldsymbol{Y},z_i) \right]\quad 
\textrm{as}\quad |u| \uparrow 1\, ,
$$
where $\Rightarrow$ denotes convergence in distribution, and
$$
\Delta(u,z)\, =\, \frac{1-\overline{u}z}{\sqrt{1-|u|^2}}\, .
$$
\end{lemma}
We may prove this result using the following simple corollary of the Lindeberg-Feller central limit theorem. 

\begin{cor}
\label{cor:LF}
Suppose that $X_0,X_1,\dots$ are i.i.d., complex valued random variables with mean zero,
satisfying condition (\ref{eq:unitvariance}).
Suppose that $\alpha_{k}(u) \in \C$ is defined for each $u \in U(0,1)$
and $k \in \{0,1,\dots \}$,
satisfying
\begin{itemize}
\item[(a)] $\frac{1}{2} \sum_{k=0}^{\infty} |\alpha_{k}(u)|^2 \to \sigma^2 > 0$ as $|u| \uparrow 1$, and
\item[(b)] $\sum_{k=0}^{\infty} |\alpha_{k}(u)|^p \to 0$ as $|u| \uparrow 1$, for some
$p>2$.
\end{itemize}
Then $\sum_{k=1}^{\infty} \textrm{Re}[\alpha_{k}(u) X_k]$ converges in distribution to 
$\sigma \chi$ as $|u| \uparrow 1$, where $\chi$ is a real standard normal  random variable.
\end{cor}
See, for example, \cite{Durrett} for the Lindeberg-Feller theorem.
With this, we can prove the lemma.

\begin{proofof}{\bf Proof of Lemma \ref{lem:clt}:}
We can write
$$
\sum_{i=1}^{n} \lambda_i \frac{f(\boldsymbol{X},\Phi(u,z_i))}{\Delta(u,z)}\,
=\, \sum_{k=0}^{\infty} \alpha_k(u) X_k\, ,
$$
where
$$
\alpha_k(u)\, =\, \sum_{i=1}^{n} \lambda_i \, \frac{\Phi(u,z_i)^k}{\Delta(u,z_i)}\, ,
$$
since $f(\boldsymbol{X},\Phi(u,z_i)) = \sum_{k=0}^{\infty} X_k \Phi(u,z_i)^k$.
A simple calculation shows that
$$
\sum_{k=0}^{\infty} |\alpha_k(u)|^2\,
=\, \sum_{i,j=1}^{n} \lambda_i \overline{\lambda}_j Q(u;z_i,z_j)\, ,
$$
where
$$
Q(u;z_i,z_j)\, =\, \frac{1}{\Delta(u,z_i) \overline{\Delta}(u,z_j)}\, \sum_{k=0}^{\infty} [\Phi(u,z_i) \overline{\Phi}(u,z_j)]^k\, .
$$
But an important property of the M\"obius transform is that it leaves the covariance of this family of random analytic functions invariant, other than multiplying by the factors $\Delta$.
In fact, this is an important property of the Gaussian analytic functions studied by Peres and Vir\'ag, since it shows that their entire distributions are stationary, as the distribution of a Gaussian process is determined by the covariance.
This is checked by summing the series to obtain
$$
Q(u;z_i,z_j)\, 
=\, \frac{1-|u|^2}{(1-\overline{u}z_i)(1-u\overline{z}_j)} \cdot \frac{1}{1-\Phi(u,z_i)\overline{\Phi}(u,z_j)}\,
=\, \frac{1}{1-z_i \overline{z}_j}\, ,
$$
which does not depend on $u$. For the same reason it shows that
$$
\textrm{var}\left(
\textrm{Re}\left[ \sum_{i=1}^{n} \lambda_i \frac{f(\boldsymbol{X},\Phi(u,z_i))}{\Delta(u,z_i)} \right]\right)\,
=\, \textrm{var}\left(\textrm{Re}\left[ \sum_{i=1}^{n} \lambda_i f(\boldsymbol{Y},z_i) \right]\right) \, ,
$$
for all $u$. This takes care of condition (a) in Corollary \ref{cor:LF}.

To check condition (b) with $p=4$, we note that
$$
\sum_{k=0}^{\infty} |\alpha_k(u)|^4\,
\leq \, n^2 \max_{i=1,\dots,n} |\lambda_i|^4 
\sum_{k=0}^{\infty} \frac{|\Phi(u,z_i)|^{4k}}{|\Delta(u,z_i)|^4}\, .
$$
But we can sum the last series for each $i$, to obtain
$$
\sum_{k=0}^{\infty} \frac{|\Phi(u,z_i)|^{4k}}{|\Delta(u,z_i)|^4}\,
=\, \frac{(1-|u|^2)^2}{|1-\overline{u}z_i|^4} \cdot \frac{1}{1-|\Phi(u,z_i)|^4}\,
=\, \frac{1-|u|^2}{(1-|z_i|^2)(|1-\overline{u}z_i|^2 + |z_i - u|^2)}\, .
$$
As long as all $z_i$ are strictly inside the unit circle, this quantity converges
to zero in the limit $|u| \uparrow 1$. This completes the proof of the lemma.
\end{proofof}

Lemma \ref{lem:clt} implies that the random analytic function
$f(\boldsymbol{X},\Phi(u,z))/\Delta(u,z)$ converges to the random analytic function $f(\boldsymbol{Y},z)$
in distribution, as $|u| \uparrow 1$, in the sense of finite dimensional marginals.
But with this we may use the following lemma of Valko and Vir\'ag from
their paper on random Schr\"odinger operators \cite{ValkoVirag}.

\begin{lemma}[Valko and Vir\'ag, 2010]
\label{lem:ValkoVirag}
Let $f_n(\omega,z)$ be a sequence of random analytic functions on a domain $D$ (which is open, connected and simply connected)
such that ${\bf E} h(|f_n(z)|) < g(z)$ for some increasing unbounded function $h$
and a locally bounded function $g$. Assume that $f_n(z) \Rightarrow f(z)$ in the sense of 
finite dimensional distributions. Then $f$ has a unique analytic version and 
$f_n \Rightarrow f$ in distribution with respect to local-uniform convergence. 
\end{lemma}

Because of this result we see that $f(\boldsymbol{X},\Phi(u,z))/\Delta(u,z)$ converges
in distribution to $f(\boldsymbol{Y},z)$ with respect to the local uniform convergence.
But by Hurwitz's theorem or Rouch\'e's theorem, this implies that the zero sets also converge in distribution,
relative to the local weak topology on point processes.
Since $\Delta(u,z)$ is finite and non-vanishing for $z \in U(0,1)$, the zero set is just the zero set
of $f(\boldsymbol{X},\Phi(u,z))$.
Combining this result with Peres and Vir\'ag's Theorem \ref{thm:PV} for the zero set of $f(\boldsymbol{Y},z)$ proves our theorem.

\section{Discussion and Extensions}

The proof presented here also may be extended to other families of Gaussian analytic functions.
In a recent book by Hough, Krishnapur, Peres and Vir\'ag \cite{HKPV}
several families of Gaussian analytic functions
were studied, whose covariances are adapted
to the classical symmetric spaces: the sphere,
the plane and the hyperbolic plane.
Some of the ensembles had been introduced before in \cite{BBL1,BBL2,DE,Kostlan,SS}.

In \cite{HKPV}, there is presented a one-parameter family of Gaussian
analytic functions, adapted to the Poincar\'e disk model of hyperbolic geometry,
for all choices of Gaussian curvature $k<0$,
as well as a model for Gaussian analytic functions on the plane corresponding to $k=0$.
They also present Gaussian polynomials adapted to the sphere for quantized values
of $k$: $k = 1/n$ for $n \in \N$.
The special case considered by Peres and Vir\'ag in \cite{PeresVirag} corresponds to
$k=-1$.
Theorem \ref{thm:main} was just for the $k=-1$ model, but extends to all the models with $k\leq 0$, with analogous proofs {\em mutatis mutandis}.
(See our original preprint \cite{LMS} for full details.)
For the spherical case,
the diameter is finite, so no limit law is possible for fixed $k=1/n$, $n \in \N$.

Finally, we end with a remark.
One could consider the joint distribution of the two analytic functions $f(\boldsymbol{X},\Phi(u_1,z))$, $f(\boldsymbol{X},\Phi(u_2,z))$ in the limit
that $u_1$ and $u_2$ both approach the boundary circle.
A simple calculation shows that for any $z_1,z_2 \in U(0,1)$,
$$
\textrm{cov}\left(\frac{f(\boldsymbol{X},\Phi(u_1,z_1))}{\Delta(u_1,z_1)}\, ,\
\frac{f(\boldsymbol{X},\Phi(u_2,z_2))}{\Delta(u_2,z_2)}\right)\\
=\, O\left( 1 + \frac{|u_1-u_2|}{\sqrt{(1-|u_1|^2)(1-|u_2|^2)}}\right)^{-1}
$$
Since the processes are Gaussian, absence of correlations
implies independence. Therefore, Lemma \ref{lem:ValkoVirag}
may be used again to conclude that the zero sets of $f(\boldsymbol{X},\Phi(u_1,z))$ and 
$f(\boldsymbol{X},\Phi(u_2,z))$ are asymptotically independent if and only if 
$\Phi(u_1,u_2) \to \infty$.

\section*{Acknowledgements}

We are very grateful to Carl Mueller and Balint Vir\'ag
for useful conversations and helpful advice.
We are also grateful to an anonymous referee for many suggestions
improving the presentation.
Parts of this research were carried out during visits of M.M.\ and S.S.\
at the University of Rochester and Memorial University of Newfoundland.
We thank these institutions for their support.

\end{document}